\documentclass[12pt]{amsart}
\usepackage{amscd,amssymb}
\usepackage[arrow,matrix]{xy}

\theoremstyle{plain}
\newtheorem{theorem}[subsection]{Theorem}
\newtheorem{lemma}[subsection]{Lemma}
\newtheorem{prop}[subsection]{Proposition}
\newtheorem{corollary}[subsection]{Corollary}

\theoremstyle{definition}
\newtheorem{remark}[subsection]{Remark}
\newtheorem{definition}[subsection]{Definition}

\numberwithin{equation}{section}
\newcommand{\A}{{\mathcal A}}
\newcommand{\I}{{\mathcal I}}
\newcommand{\B}{{\mathcal B}}
\newcommand{\M}{{\mathcal M}}
\newcommand{\E}{{\mathcal E}}

\newcommand{\Z}{\mathbb{Z}}
\newcommand{\Q}{\mathbb{Q}}

\newcommand{\C}{\mathbb{C}}

\newcommand{\K}{\mathbb{K}}

\newcommand{\PP}{\mathbb{P}}

\DeclareMathOperator{\Hom}{Hom}

\DeclareMathOperator{\spn}{span}
\DeclareMathOperator{\id}{id}

\DeclareMathOperator{\der}{Der}
\DeclareMathOperator{\aut}{Aut}

\newcommand{\surj}{\twoheadrightarrow}

\newenvironment{romenum}
{

\begin{enumerate}}{\end{enumerate}}

\begin{document}

\title[Geodesics and derivations]{Isometry-invariant geodesics and\\
nonpositive derivations of the cohomology}

\author[\c{S}tefan Papadima]{\c{S}tefan Papadima$^{* \dagger}$}

\author[Lauren\c{t}iu P\u{a}unescu]{Lauren\c{t}iu P\u{a}unescu$^\dagger$}

\thanks{$^*$Partially supported by grant 
CERES 152/2003 of the Romanian Ministry of
Education and Research.}

\thanks{$^\dagger$Partially supported by grant U4249 Sesqui R\&D/2003
of the University of Sydney.} 



\begin{abstract}
 We introduce a new class of zero-dimensional weighted complete intersections,
 by abstracting the essential features of $\Q$-cohomology algebras of equal 
 rank
 homogeneous spaces of compact connected Lie groups. We prove that, on a 
 1-connected closed manifold $M$ with $H^*(M,\Q)$ belonging to this class,
 every isometry has a non-trivial invariant geodesic, for any metric on $M$.
 We use $\Q$-surgery to construct large classes of new examples for which 
 the above result may be applied.
\end{abstract}

\maketitle

\section{Introduction}
\label{sec:intro}

\subsection{Complexes which look like homogeneous spaces.}
\label{subsec=intro1}
Let $\A$ be a {\it weighted, zero--dimensional} (that is, artinian), 
{\it complete intersection}
($WACI$), i.e., a commutative graded $\Q$-algebra of the form
\begin{equation}
\label{eq:WACI}
\A = \Q[x_1, \dots, x_n]/ \I \, ,
\end{equation}
where the variables $x_i$ have positive even weights, 
$w_i:= \mid x_i \mid$, 
and the ideal $\I$ is generated by a regular sequence,
\begin{equation}
\label{eq:reg}
\I = (f_1, \dots, f_n) \, ,
\end{equation}
of weighted-homogeneous polynomials, $f_i$.

We are going to introduce a special class of $WACI$'s, 
with an eye for applications in Riemannian geometry.

One may speak about the {\it pseudo-homotopy groups} of a
$1$-connected, finitely generated, graded-commutative $\Q$-algebra, $\A$:
\begin{equation}
\label{eq:pseudo}
\pi^*_*(\A) = \bigoplus_{i\ge 0,j>1} \pi^j_i(\A) \, .
\end{equation}

These algebraic invariants of $\A$ are finite-dimensional 
graded $\Q$-vector spaces, $\pi_i(\A) := \oplus_j \pi^j_i(\A)$, 
for all $i\ge 0$. See \cite{HS}.

When $\A$ is a $WACI$, as above, one knows that $\pi_{>1}(\A) =0$.
Moreover, for any $1$-connected $CW$-complex, $S$, such that
$H^*(S, \Q)=\A^*$, as graded algebras, one has the following 
topological interpretation (see Sullivan~\cite{S}):
\begin{equation}
\label{eq:pi}
\begin{cases}
\pi^*_0(\A) = \Hom_{\Q} (\pi_{* = \text{even}}(S)\otimes \Q, \Q), 
& \text{and}\\
\pi^*_1(\A) = \Hom_{\Q} (\pi_{* = \text{odd}}(S)\otimes \Q, \Q) & .
\end{cases}
\end{equation}

Assuming $\A \neq \Q$, that is, $\pi_1^*(\A) \neq 0$, set now:
\begin{equation}
\label{eq:deftop}
k_{\A} := \text{max} \, \{ k \mid \pi^{2k-1}_1(\A) \ne 0 \} \, .
\end{equation}

Denote by $\der^* (\A)$ the graded Lie algebra of 
homogeneous derivations of $\A^*$, with Lie bracket 
given by graded commutator, and note that $\der^*(\Q)=0$.

\begin{definition}
\label{def:simple}

We shall say that the $WACI$ $\A$ is {\em simple} if:
\begin{romenum}
\item \label{derneg}
$\der^{<0}(\A) =0$.
\item \label{derzero}
$\dim_{\Q} \der^0(\A) =1$.
\item \label{top}
$\dim_{\Q} \pi^{2k_{\A}-1}_1(\A) =1$.
\end{romenum}
\end{definition}

\begin{definition}
\label{def:hh}

A 1-connected $CW$-complex $S$ is said to be {\em homologically
homogeneous} ($HH$) if 
$$ H^*(S, \Q) = \bigotimes_{j\in J} \A^*_j \, ,$$
as graded algebras, where the index set $J$ is finite and 
the algebra $\A_j$ is a simple $WACI$, for every $j\in J$.
\end{definition}

Both definitions are motivated by geometry. Firstly, it is
well-known that the $\Q$-cohomology algebra of an equal rank
homogeneous space, $G/K$, of compact connected Lie groups,
is a $WACI$; see Borel~\cite{B}. Moreover, $G/K = \prod_{j \in J}
G_j/K_j$, where in addition each Lie group $G_j$ is simple.
Now, for $G$ simple and $\A^* = H^*(G/K, \Q)$, properties
\eqref{derneg} and \eqref{derzero} from Definition~\ref{def:simple}
were proved by Shiga and Tezuka in \cite{ST}. 
Property \eqref{top} follows from
\cite[Theorem 1.3]{P2}, via \eqref{eq:pi}. Therefore, all equal
rank homogeneous spaces $G/K$ mentioned above are $HH$, in the 
sense of our Definition~\ref{def:hh}.

We exhibit in Section~\ref{sec=new} 
(Theorems~\ref{thm:smooth} and~\ref{thm:smooth2}) families of
new examples of $HH$ closed manifolds. Their $\Q$--cohomology algebras
are actually simple $WACI$'s, not isomorphic to the $\Q$--cohomology algebra
of any equal rank homogeneous space $G/K$, with $G$ not necessarily
simple. We first use tools from our previous work~\cite{PP}, to obtain
simple $WACI$'s from $1$-dimensional {\em reduced} weighted 
complete intersections. Rational surgery methods, due to Sullivan~\cite{S}
and Barge~\cite{Ba}, enable us to realize them by closed manifolds.

\subsection{Isometry--invariant geodesics.}
\label{subsec=intro2}

Let $M^m$ be a 1-connected, closed, Riemannian $m$-manifold.
A basic problem is to find topological conditions
on $M$ guaranteeing that

\begin{equation}
\label{eq:ideal}
\text{every isometry $a$ has a non-trivial $a$-invariant
geodesic.}
\end{equation}

Recall that a geodesic curve $\gamma$ is called $a$-invariant 
if there exists a period $t$ such that
$a(\gamma (x))= \gamma (x+t)$, for any $x$. When $a =\id$, one 
recovers the classical notion of closed geodesic. Strong existence results
for closed geodesics were obtained, under simple topological
hypotheses on $M$, by Sullivan and Vigu\' e--Poirrier in~\cite{SV},
using rational homotopy methods.

Further refinements led to
various answers to the general fundamental question \eqref{eq:ideal}. Recall
that a $1$-connected finite complex $S$ is called {\it elliptic} 
if $\dim_{\Q} \pi_*(S)\otimes \Q <\infty$, and {\it hyperbolic}
otherwise. Grove and Halperin proved in \cite{GH} that property
\eqref{eq:ideal} holds for all hyperbolic $M$, and for all
elliptic $M$ with $m$ odd. 

On the other hand,  $HH$ manifolds are elliptic
(see \eqref{eq:pi}) and even-dimensional (see Halperin~\cite{H}).
In the particular case of equal rank homogeneous spaces,
$M=G/K$,  property \eqref{eq:ideal} holds, when $K$ is
a maximal torus (\cite[Theorem 1.4(i)]{P1}) or $G$ is
simple (\cite[Corollary 1.4(i)]{P2}).

Our main result below significantly enlarges the class
of manifolds $M$ having the ideal property \eqref{eq:ideal}.

\begin{theorem}
\label{thm:main}
If the closed manifold $M$ is homologically homogeneous
(in the sense of {\em Definition \ref{def:hh}}), then, 
given an arbitrary metric on $M$, there is
a non-trivial $a$-invariant geodesic, for every isometry, $a$.
\end{theorem}

Along the way, we shall also clarify the structure of the
automorphism group of certain related graded algebras, $\B^*$.
More precisely, assume that $\B^* =\otimes_{j\in J} \B_j^*$,
where $J$ is finite, and each connected, commutative, evenly--graded
$\C$--algebra $\B^*_j$ is artinian and satisfies properties 
\eqref{derneg} and \eqref{derzero} from 
Definition~\ref{def:simple} (over $\C$). Then the connected 
component of $1$ in the (linear algebraic) group of graded algebra
automorphisms of $\B^*$ is an algebraic torus, $(\C^*)^r$, where
$r= \mid J\mid$. See Corollary~\ref{cor=struct1}.

\section{Proof of the main result}
\label{sec=pfmain}

If $a$ has no invariant geodesic, then both
$\pi_*(a)\otimes \C$ and $\id -\pi_*(a)\otimes \C$ 
are unimodular elements of
$GL (\pi_*(M)\otimes \C)$; see \cite{G}. Set
$k_j =k_{\A_j}$, for $j\in J$, and 
$k = \text{max} \, \{ k_j \mid j\in J \}$.
We will derive a contradiction, for $*= 2k -1$.

\subsection{Dual homotopy representations.}
\label{subsec=step1}
Rational homotopy theory methods~\cite{S} may be used to get 
information on $\pi_{2k-1}(a)\otimes \C$. Tools from 
the representation theory of linear algebraic groups
(see e.g. \cite{Hu}) will provide additional information, 
so we will use $\C$-coefficients.

Set $\A^*= H^*(M, \Q)$. Since $M$ is $HH$, $\A$ is a $WACI$.
As shown in \cite{S}, this implies that the $\C$-minimal model
of $M$, $(\M, d)$, has the following properties. It is a bigraded
algebra of the form $\M = \C [Z_0]\otimes \bigwedge Z_1$, where
$Z_0= \C-\spn \{ x_1,\dots , x_n \}$, with $x_i$ of even upper degree,
and $Z_1= \C-\spn \{ y_1, \dots ,y_n \}$, with $y_i$ of odd upper degree.
The differential $d$ is homogeneous, with upper degree $+1$ and
lower degree $-1$. More precisely, $dx_i=0$ and
$dy_i \in \Q [x_1, \dots, x_n]$ is a polynomial with no linear part,
for all $i$. 
Moreover, $H_+ (\M ,d)=0$, and $Z_j^*= \pi_j^*(\A) \otimes \C$;
see \cite{HS}. 

Denote by $\E$ the group of self-homotopy equivalences of $M$. By
\eqref{eq:pi}, the natural action of $\E$ on $\pi_{2k-1}(M)$ gives rise 
to a dual homotopy (anti)representation, 
\begin{equation}
\label{eq=pirep}
\tau \colon \E \longrightarrow GL (Z_1^{2k-1})\, .
\end{equation}

Denote by $\I \subset \C [Z_0]$ the ideal generated by 
$\{ dy_i \}_{1\le i\le n }$, and by $\C^+ [Z_0]$ the strictly positively
graded part of $\C [Z_0]$.

\begin{lemma}
\label{lem=step11}

The differential $d$ induces a vector space isomorphism,
\[
\Psi \colon Z_1^{2k-1}\stackrel{\sim}{\longrightarrow}
[\I/ \C^+ [Z_0]\cdot \I]^{2k}\, .
\]
\end{lemma}

\begin{proof}
The surjectivity of $\Psi$ is immediate. To check injectivity, assume that
\[
\sum^n_{i=1} c_i dy_i =\sum^n_{i=1} p_i dy_i\, ,
\]
with $c_i\in \C$ and
$p_i \in \C^+ [Z_0]$. Since $H_1(\M)=0$, we infer that
$\sum^n_{i=1} (c_i-p_i) y_i \in d(\M_2)$. Plainly, 
$d(\M_2)\subset \C^+ [Z_0]\otimes Z_1$, hence $c_i=0$, for all $i$.
\end{proof}

Denote by $\aut_{\C}(\A)$ the group of graded algebra automorphisms
of $\A^*\otimes \C$. It is a linear algebraic group, defined over $\Q$.
There is a cohomology (anti)representation,
\begin{equation}
\label{eq=hrep}
h \colon \E \longrightarrow \aut_{\C}(\A)\, ,
\end{equation}
which associates to $a\in \E$ the cohomology automorphism
$H^*(a, \C)$.

Our next goal is to relate \eqref{eq=pirep} and \eqref{eq=hrep}, by
constructing an algebraic representation,
\begin{equation}
\label{eq=psrep}
\nu \colon \aut_{\C} (\A)\longrightarrow GL (Z_1^{2k-1})\, ,
\end{equation}
called the {\em pseudo--homotopy representation}. To do this, we first
consider the linear algebraic group $\aut_{\I}(\C [Z_0])$, consisting
of those graded algebra automorphisms, $\alpha$, of $\C [Z_0]$, which
leave the ideal $\I$ invariant. From Lemma~\ref{lem=step11}, we get
a natural algebraic representation,
\begin{equation}
\label{eq=irep}
\nu \colon \aut_{\I}(\C [Z_0])\longrightarrow GL (Z_1^{2k-1})\, .
\end{equation}

Note that $\A^*\otimes \C= H^*(\M ,d)= H_0^*(\M , d)= \C [Z_0]/\I$.
It follows that $\aut_{\C}(\A)$ is the quotient of 
$\aut_{\I}(\C [Z_0])$ by the subgroup of those elements $\alpha$
having the property that
\begin{equation}
\label{eq=alfa}
\alpha (x_i) \equiv x_i \quad \text{(modulo $\I$)}\, , \quad
\text{for all} \quad i\, .
\end{equation}

It is easy to infer from \eqref{eq=alfa} that $\alpha$ induces 
the identity on $\I/\C^+ [Z_0] \cdot \I$, by resorting to the minimality
property of $(\M ,d)$. Thus, \eqref{eq=irep} factors to give 
the desired representation, \eqref{eq=psrep}.

\begin{lemma}
\label{lem=step12}

The representations \eqref{eq=pirep}, \eqref{eq=hrep} and
\eqref{eq=psrep} are related by:
\[
\tau \quad = \quad \nu \circ h \, \, .
\]
\end{lemma}

\begin{proof}
For $a\in \E$, denote by $\widehat{a} \colon (\M ,d) \to (\M ,d)$
the Sullivan minimal model of $a$. It is a differential graded algebra
($DGA$) automorphism, with the property that $H^*(\widehat{a})= h(a)$.

Recall from \cite{S} that $\tau (a)=\pi^{2k-1}(\widehat{a})$, where
$\pi^{2k-1}(\widehat{a}) \colon 
Z_1^{2k-1} \to Z_1^{2k-1}$ denotes the map induced by
the restriction of $\widehat{a}$ to $Z_1^{2k-1}$, modulo 
decomposable elements of $\M^+$. 

To compute $\nu(h(a))$, we need to lift $H^*(\widehat{a})$ to a
graded algebra map, $\C [Z_0] \to \C [Z_0]$, leaving $\I$
invariant. To do this, we may proceed as follows. Start by noting
that $\widehat{a}$ increases lower degrees. It is enough to check
this on algebra generators of $\M$, which is obvious on $Z_0$; on
$Z_1$, this follows from the fact that $\widehat{a}(y_i)$ has odd upper
degree, for all $i$. Write then 
$\widehat{a}= \sum_{j\ge 0} \widehat{a}_j$,
with $\widehat{a}_j \in \Hom_{\C}(\M_*, \M_{*+j})$. An easy lower
degree argument, together with the bihomogeneity property of $d$,
show that $\widehat{a}_0$ is a bigraded differential algebra map.
Set $\alpha :=\widehat{a}_0 \colon \C [Z_0] \to \C [Z_0]$. Since
$H_+(\M ,d)=0$, it is easily seen that $\alpha \in \aut_{\I}(\C [Z_0])$
lifts $H^*(\widehat{a})\in \aut_{\C}(\A)$. 

To finish our proof, let us check that
\begin{equation}
\label{eq=step12}
d \pi^{2k-1}(\widehat{a})(y_i)\equiv \widehat{a}_0 (dy_i)\, ,
\quad \text{modulo}\quad \C^+[Z_0]\cdot \I \, ,
\end{equation}
for $y_i$ of (upper) degree $2k-1$; see Lemma~\ref{lem=step11}.
Note that$\pi^{2k-1}(\widehat{a})(y_i)$ has odd upper degree; hence,
$\pi^{2k-1}(\widehat{a})(y_i)\equiv \widehat{a}_0(y_i)$, modulo
$\C^+ [Z_0]\otimes Z_1$. Applying $d$, we get \eqref{eq=step12}.
\end{proof}

\subsection{Derivations and the structure of $\aut_{\C}(\A)$.}
\label{subsec=step2}
Denote by $\aut_{\C}^1(\A)$ the connected component of
$1\in \aut_{\C}(\A)$. We are going to exploit conditions
\eqref{derneg} and \eqref{derzero} from Definition~\ref{def:simple},
to give a complete description of the algebraic group
$\aut_{\C}^1(\A)$. 

In this subsection, all graded-commutative algebras $\B^*$ will be
supposed to be connected, finite-dimensional over a field $\K$, 
and evenly-graded. The examples we have in mind are $\B=\A\otimes \C$, 
where $\A$ is a $WACI$. For $p\in \Z$, recall that 
\[
\der^p(\B) := \{ \theta\in \Hom_{\K}(\B^*, \B^{*+p})\, \mid \, 
\theta (ab) =\theta a\cdot b + a\cdot \theta b\, , \forall a,b\in \B \}\, .
\]

We will be particularly interested in the Lie algebra $\der^0(\B)$.
In our basic examples, one knows that $\der^0(\A \otimes \C)$ is the
Lie algebra of $\aut_{\C}^1(\A)$; see~\cite[12.5 and 13.2]{Hu}. 
We first explore the consequences of the vanishing of $\der^{<0}$
on the structure of $\der^0$.

\begin{lemma}
\label{lem=dero}

Assume that $\der^{<0}(\B_j)=0$, for $j=1,\dots ,r$. Then 
the Lie algebra $\der^0 (\otimes_{j=1}^r \B_j)$ is isomorphic to
the direct product $\prod_{j=1}^r \der^0(\B_j)$.
\end{lemma}

\begin{proof}
Consider $\theta\in \der^p(\B_1\otimes \B_2)$, with $\B_1$, $\B_2$ and $p$
arbitrary. Pick $\{ b_{\beta} \}_{\beta}$, a homogeneous $\K$-basis
of $\B_2$, comprising $1$. Write $\theta_{\mid \B_1\otimes 1} =
\sum_{\beta} \theta_{\beta}\otimes b_{\beta}$. A straightforward
computation shows that 
$\theta_{\beta}\in \der^{p-\mid b_{\beta}\mid}(\B_1)$,
for every $\beta$. This shows inductively that 
$\der^{<0}(\otimes_{j=1}^r \B_j)=0$, if $\der^{<0}(\B_j)=0$, for all $j$.

Let again $\B_1$ and $\B_2$ be arbitrary. It is easy to check that the map
\begin{equation}
\label{eq=deco}
\Phi \colon \der^0(\B_1) \times \der^0(\B_2)\longrightarrow
\der^0(\B_1\otimes \B_2)\, ,
\end{equation}
defined by sending $(\theta_1, \theta_2)$ to
$\theta_1\otimes \id + \id\otimes \theta_2$, is a Lie algebra
monomorphism. We claim that $\Phi$ is an isomorphism, whenever
$\der^{<0}(\B_1)=\der^{<0}(\B_2)=0$. To check surjectivity, pick
$\theta\in \der^0(\B_1\otimes \B_2)$ arbitrary. 
By the discussion from the preceding paragraph, 
$\theta_{\mid \B_1\otimes 1}=\theta_1\otimes 1$, 
and $\theta_{\mid 1\otimes \B_2}=1\otimes \theta_2$, with
$\theta_j\in \der^0(\B_j)$, for $j=1,2$. We infer that $\theta =
\theta_1\otimes \id + \id\otimes \theta_2$, which proves our claim. By
making repeated use of the isomorphism \eqref{eq=deco}, one verifies
inductively the asertion of the Lemma.
\end{proof}

Let now $\K$ be algebraically closed, of characteristic zero.
For arbitrary $\B_1, \dots ,\B_r$, one has an algebraic group morphism,
with source the algebraic $r$--torus, $(\K^*)^r$, 
\begin{equation}
\label{eq=torus}
\rho \colon (\K^*)^r \longrightarrow \aut^1(\bigotimes_{j=1}^r \B_j)\, ,
\end{equation}
defined by $\rho (t)=\otimes_{j=1}^r \rho_j (t_j)$, for
$t= (t_1, \dots ,t_r)$, where each {\em grading automorphism}
$\rho_j(t_j)$ acts on $\B_j^q$ as $t_j^q \cdot \id$.

Our next result, which extends Corollary $(ii)$ from~\cite{PP}, 
will provide an important step in the proof of Theorem~\ref{thm:main}.

\begin{corollary}
\label{cor=struct1}

Assume that each $\B_j$ satisfies properties \eqref{derneg} and
\eqref{derzero} from {\em Definition~\ref{def:simple}} (over $\K$).
Then the above morphism $\rho$ is onto, with finite kernel. In 
particular, $\aut^1(\otimes_{j=1}^r \B_j)$ is an $r$--torus.
\end{corollary}

\begin{proof}
Property \eqref{derzero} forces $\B_j \neq \K$, for all $j$. This
readily implies that the kernel of $\rho$ is finite. Together with
Lemma~\ref{lem=dero}, this in turn guarantees the surjectivity of
$\rho$, given that $\der^0(\otimes_{j=1}^r \B_j)$ is the Lie algebra 
of $\aut^1(\otimes_{j=1}^r \B_j)$~\cite{Hu}, via a dimension argument.
The last assertion of the Corollary follows from \cite[16.1--2]{Hu}.
\end{proof}

\subsection{A convenient $\Q$-basis.}
\label{subsec=step3}
To prove Theorem~\ref{thm:main}, we have seen, at the beginning 
of this section, that it is enough to show that one cannot 
have simultaneously $\det \tau(a)=\pm 1$ and $\det(\id -\tau(a))=\pm 1$,
where $a\in \E$ is an arbitrary self-homotopy equivalence of $M$. 
To this end, the major step consists in finding a $\Q$-basis of
$\pi_1^{2k-1}(\A)$, the $\Q$-form of $Z_1^{2k-1}$, with respect to which
$\tau (a)$ has a particularly simple matrix, for {\em every} $a\in \E$.
This will be done by resorting to the last property from 
Definition~\ref{def:simple}, and then using Corollary~\ref{cor=struct1}.

Set  $J' := \{ j\in J \mid k_j=k \}$. One knows that 
$\pi_1^{2k-1}(\A)= \oplus_{j\in J} \pi_1^{2k-1}(\A_j)$; 
see \cite{HS}, \cite{S}. With the notation from \S~\ref{subsec=step1},
pick $y_j \in \pi_1^{2k-1}(\A_j)$, $y_j \neq 0$, for $j\in J'$. Then
obviously $\{ y_j \}_{j\in J'}$ will be a basis of the $\Q$-form
of $Z_1^{2k-1}$. 

\begin{lemma}
\label{lem=qmat}

For any $a\in \E$, there exist $\{ \lambda_j\in \Q \}_{j\in J'}$,
and a permutation $\sigma$ of $J'$, such that
\[
\tau (a) y_j = \lambda_j y_{\sigma (j)}\, ,
\quad \text{for} \quad j\in J' \, .
\]
\end{lemma}

\begin{proof}
By Lemma~\ref{lem=step12}, it suffices to prove the above statement,
replacing $\tau (a)$ with $\nu (\alpha)$, where 
$\alpha\in \aut_{\C}(\A)$ is arbitrary, and without demanding the
rationality of the $\lambda$-vector. (For $\alpha= h(a)$, the rationality
property easily follows, since $\tau (a)$ respects by construction 
the $\Q$-structure.)

Since $\aut_{\C}(\A)$ obviously normalizes $\aut_{\C}^1(\A)$, the map
on characters induced by conjugation by $\alpha^{-1}$ permutes 
the weights of the restriction of the linear representation 
$\nu$ to  $\aut_{\C}^1(\A)$, and $\nu (\alpha)$ permutes the corresponding
weight spaces; see~\cite[11.4]{Hu}. Due to the surjectivity of $\rho$
from \eqref{eq=torus}, these weights coincide with the weights of
$\nu \rho$, and the corresponding weight spaces are equal; 
see again~\cite[11.4]{Hu}.

At the same time, the structure of the representation $\nu \rho$
is easy to describe. Indeed, we claim that, for $j\in J'$, $y_j$ is
a weight vector of $\nu \rho$, with weight $t_j^{2k}$. This 
may be checked without difficulty, starting from the fact that
$(\M ,d)$, the minimal model of the $DGA$ $(\A \otimes \C, d=0)$,
splits as $\otimes_{j\in J} (\M_j ,d_j)$, where $(\M_j ,d_j)$ is
the minimal model of $(\A_j \otimes \C ,d=0)$, according to~\cite{S}.

We infer that the (distinct) weights of $\nu \rho$ are
$\{ t_j^{2k} \}_{j\in J'}$, with one--dimensional corresponding
weight spaces, spanned by $\{ y_j \}_{j\in J'}$. The desired form of 
$\nu (\alpha)$
is thus obtained.
\end{proof}

\subsection{End of proof of Theorem~\ref{thm:main}.}
\label{subsec=step4}
We may now compute the characteristic polynomial of $\tau (a)$, 
$P(z)$, as follows. Set $s= \mid J' \mid$, let 
$\sigma_1, \dots ,\sigma_m$
be the associated cycles of $\sigma$, with lengths $s_1, \dots ,s_m$,
and put $\gamma_i := \prod_{j\in \text{supp} (\sigma_i)} \lambda_j$,
for $i=1, \dots ,m$. Then:
\begin{equation}
\label{eq=charp}
P(z)\quad = \quad \prod_{i=1}^m (z^{s_i}-\gamma_i)\,\, .
\end{equation}

We claim that
\begin{equation}
\label{eq=fclaim}
\gamma_i \in \Z\,\, ,\quad \text{for all} \quad i \,\, .
\end{equation}

Granting the claim for the moment, we may quickly finish the proof
of Theorem~\ref{thm:main}. From the assumption that the isometry 
$a$ has no non-trivial invariant geodesic, we have deduced that both
$P(0)$ and $P(1)$ are equal to $\pm 1$. Due to \eqref{eq=fclaim},
this is impossible, and we are done.

Coming back to the claim, let us remark that the polynomial $P(z)$
from \eqref{eq=charp} actually belongs to $\Z [z]$, since $\tau (a)$
also preserves the natural $\Z$-structure of $\pi_1^{2k-1}(\A)$; see
\eqref{eq:pi}. The following elementary lemma will thus complete our proof.

\begin{lemma}
\label{lem=zmat}

Let $P(z)$ be given by \eqref{eq=charp}.
If $P(z)\in \Z [z]$, then $\gamma_i\in \Z$, for all $i$.
\end{lemma}

\begin{proof}
Write $\gamma_i =\frac{u_i}{v_i}$, with $u_i$ and $v_i$ relatively prime.
We have
\[
v_1 \cdots v_m P(z) = \prod_{i=1}^m (v_i z^{s_i}- u_i) \, .
\]
Applying the Gauss Lemma (see e.g. \cite[p.127]{L}) in the above equality,
we infer that $v_1 \cdots v_m =\pm 1$, since the content of 
the monic polynomial $P(z)$ is $1$.
\end{proof}

\section{Reduced complete intersection and smoothing}
\label{sec:smooth}

To construct new examples of application of Theorem~\ref{thm:main},
we will proceed in two steps. We will first present a general way of
obtaining simple $WACI$'s. Secondly, we will review the obstruction theory
associated to the problem of realizing Artinian graded $\Q$-algebras
by smooth manifolds.

\subsection{Reduced complete intersection.}
\label{subsec:red}
The construction of simple $WACI$'s 
is based on the Corollary on 
page 597 of our previous work \cite{PP}.

\begin{prop}
\label{lem:red}

Let $\A$ be a $WACI$, as in {\em \eqref{eq:WACI}--\eqref{eq:reg}},
where the $f$-weights satisfy 
$0< \mid f_1 \mid\le \dots \le \mid f_{n-1} \mid <\mid f_n \mid$,
and the polynomial $f_n$ has no linear part. If 
$\Q [x_1, \dots ,x_n]/(f_1, \dots ,f_{n-1})$ is reduced
(over $\C$), then $\A$ is simple.
\end{prop}

\begin{proof}
Conditions \eqref{derneg} and \eqref{derzero} follow
directly from the abovementioned Corollary, and its proof;
see \cite{PP}.

To check the last condition \eqref{top}, one has just to
recall from \cite{S} the recipe for computing $\pi^*_1$
of an arbitrary $WACI$. Set $X :=\Q -\spn \, \{x_1, \dots, x_n \}$
and $Y := \Q -\spn \, \{y_1, \dots, y_n \}$, with 
$\mid y_i \mid = \mid f_i \mid -1$. Then:
\begin{equation}
\label{eq:min}
\pi^*_1(\A) = \ker \, (L) \, ,
\end{equation}
where $L : Y\to X$ sends $y_i$ to the linear part of $f_i$.
\end{proof}

\subsection{The smoothing problem.}
\label{subsec=smpb}
Let $\A^*$ be a $1$-connected {\em Poincar\' e duality algebra} ($PDA$)
over $\Q$, of formal dimension $m$. That is, a finite-dimensional
graded-commutative $\Q$-algebra with $\A^1=0$, $\A^m\cong \Q$ and
$\A^{>m}=0$, for which the multiplication gives a duality pairing,
$\A^i\otimes \A^j \to \Q$, for $i+j=m$. We shall say that $\A$ is
{\em smoothable} if $\A^*\cong H^*(M^m ,\Q)$, as graded $\Q$-algebras,
where $M$ is a $1$-connected closed smooth $m$-manifold. A natural
source of $1$-connected $PDA$'s of even formal dimension
is provided by $WACI$'s; see~\cite[Theorem 3]{H}.

Let $\A$ be an arbitrary $1$-connected $PDA$, of formal dimension $m$.
If $\A$ is smoothable, we may pick an orientation class,
$[M]\in H_m(M ,\Z)$, which defines an algebraic orientation,
$\omega\in \A^m \setminus \{ 0 \}$. Likewise, the Pontrjagin classes
$p_i(M)$ define a total algebraic Pontrjagin class,
$q =1+ \sum_{i\ge 1} q_i$, with $q_i\in \A^{4i}$. If $m=4k$,
$\omega$ and $q$ may be used to describe the following three
well-known obstructions to smoothability (see \cite{M}).

By construction, $\{ \langle q_1^{e_1}\cdots q_k^{e_k}\, , \, 
\omega \rangle \, \mid \, \sum_{i=1}^k ie_i=k \}$ are the
{\em Pontrjagin numbers} of a closed oriented smooth $4k$-manifold.

Secondly, the nondegenerate quadratic form $\varphi (\A ,\omega)$, 
defined on $\A^{2k}$ by Poincar\' e duality and the orientation,
must be a sum of signed squares, over $\Q$ (the {\em integrality
property}).

Finally, the {\em signature formula} must hold:
\begin{equation}
\label{eq=sign}
\sigma (\varphi (\A , \omega))\quad =\quad \langle L_k (q_1,\dots ,q_k)
\, , \,\omega \rangle \, ,
\end{equation}
where $\sigma$ denotes the signature and $L_k$ is the $k$-th
Hirzebruch polynomial.

It turns out that these obstructions completely control 
the smoothing problem, as explained in the proposition below, 
which is an immediate consequence
of a basic $\Q$-surgery result due to D. Sullivan~\cite[Theorem 13.2]{S}
(see also~\cite[Theorem 8.2.2]{Ba}). 

\begin{prop}
\label{prop=surger}
Let $\A$ be a $1$-connected $\Q$--$PDA$, of formal dimension $m$.
If $m\neq 4k$, then $\A$ is smoothable. If $m=4k$, then $\A$ is
smoothable if and only if there exist 
$\omega\in \A^{4k}\setminus \{ 0 \}$
and $\sum_{i\ge 1} q_i \in \oplus_{i\ge 1} \A^{4i}$, which verify 
the abovementioned properties, concerning Pontrjagin numbers,
integrality, and the signature formula.
\end{prop}

\begin{proof}
Realize $\A$ by a $1$-connected, $\Q$-local, Poincar\' e complex $S$.
Everything follows then by applying the results quoted from \cite{S}
and \cite{Ba} to $S$, except the fact that, for $m=4$, 
smoothability holds, as soon as the obstructions are satisfied. In
this very easy case, the Poincar\' e quadratic form on $\A^2$ is a sum
of $r$ squares minus a sum of $s$ squares, by integrality. Then plainly
$\A^*= H^*((\#_r \C \PP^2)\# (\#_s \overline{\C \PP^2}), \Q)$
(where $\overline{\C \PP^2}$ denotes $\C \PP^2$ with the 
opposite complex orientation), for $r+s> 0$, and $\A^*=H^*(S^4 ,\Q)$,
for $r+s=0$.
\end{proof}

\section{$HH$ manifolds which are not homogeneous}
\label{sec=new}

According to the signature of their cohomology algebras, our new
$HH$ manifolds fall into two classes. We begin with the (easier)
zero signature case.

\subsection{Split simple algebras.}
\label{subsec=split}
There is a large freedom of choice of (discrete) parameters
for our construction. More precisely:
pick any $n\ge 2$, $k\ge 1$, then arbitrary positive even weights,
$\{ w_i \}_{1\le i\le n}$, and integers $a_i\ge 2$ such that
$w_i a_i=d$, for $i=1, \dots, n$. Set:
\begin{equation}
\label{eq:f}
\begin{cases}
f_1 = x^{a_1}_1 - x^{a_2}_2 \\
\dots  \\
f_{n-1} = x^{a_{n-1}}_{n-1} - x^{a_n}_n \\
f_n = x^{2k  a_n}_n 
\end{cases}
\end{equation}

\begin{prop}
\label{prop:discr}

Define $\A := \Q [x_1, \dots, x_n]/ (f_1, \dots, f_n)$,
where $f_1, \dots, f_n$ are given by \eqref{eq:f} above. Then:
\begin{enumerate}
\item \label{part1}
$\A$ is a simple $WACI$.
\item \label{part2}
$\A$ is a $1$-connected Poincar\'e duality algebra,
with even formal dimension.
\item \label{part3}
The non-trivial odd pseudo-homotopy groups of $\A$ are:
\[
\begin{cases}
\pi^{d-1}_1(\A) = \Q^{n-1} \, , & {\rm and}\\
\pi^{2kd -1}_1(\A) = \Q & .
\end{cases}
\]
\end{enumerate}
\end{prop}

\begin{proof}
Part \eqref{part1} follows from Proposition~\ref{lem:red}, via
Lemma 2.1(ii) from \cite{PP}. Part \eqref{part3} is a
direct consequence of \eqref{eq:min}. For Part \eqref{part2}, see
\cite[Theorem 3]{H}.
\end{proof}

Smoothability will follow from the next well-known lemma, whose proof
is included for the reader's convenience.

\begin{lemma}
\label{lem:pzero}

Let $\A$ be a $1$-connected, commutative, Poincar\' e
duality algebra over $\Q$, with even formal dimension.
Assume that there is a $\Q$-subspace, $N\subset \A$,
with $\dim_{\Q} \A = 2 \dim_{\Q}N$, and such that
$N \cdot N=0$, where the dot denotes the Poincar\' e 
inner product. Then there is a $1$-connected, closed, smooth
manifold $M$, with zero signature, such that 
$\A^* =H^*(M, \Q)$, as graded algebras.
\end{lemma}

\begin{proof}
Immediate, using Proposition~\ref{prop=surger}. When the formal
dimension of $\A$ is $4k$, we may use any orientation,
$\omega \in \A^{4k}\setminus \{ 0 \}$, and  the trivial 
Pontrjagin class, $q=1$. Now, Poincar\' e duality implies that the
inner product space $(\A, \cdot)$ equals
$(\A^{2k}, \cdot)$, in the Witt group $W(\Q)$. On the other
hand, our assumptions on $N$ mean precisely that $(\A, \cdot)$
equals zero in $W(\Q)$. See \cite[I.6--7]{MH}.
In particular, $\sigma(\A^{2k}, \cdot)=0$.
With these remarks, it is now easy to check 
the three obstructions to smoothing.
\end{proof}

\begin{theorem}
\label{thm:smooth}

Let $\A$ be a simple $WACI$ as in 
{\em Proposition \ref{prop:discr}}. Then:
\begin{enumerate}
\item \label{mfd}
$\A^*$ is the $\Q$-cohomology algebra of a  
$1$-connected, closed smooth manifold, $M$.
The signature of $M$ is zero.
\item \label{nothom}
If $n\ge 4$, $\A^*$ is not isomorphic to the
cohomology algebra of any equal rank homogeneous space,
$G/K$.
\end{enumerate}
\end{theorem}

\begin{proof}
Part \eqref{mfd}. According to Lemma
\ref{lem:pzero}, we are done, as soon as we find 
an isotropic subspace, $N\subset \A$, having the appropriate
dimension. This in turn may be done in the following way.

Set $\B = \Q[y_1, \dots, y_n]/(y_1-y_2, \dots, y_{n-1}-y_n, y^{2k}_n)$,
where $\mid y_i \mid =d$, for all $i$. Plainly, 
$\B = \Q[y_0]/(y^{2k}_0)$, with $\mid y_0 \mid =d$.
Consider the graded algebra map, $\phi :\B^* \to \A^*$,
which sends  $y_i$ to $x^{a_i}_i$, for $i=1, \dots, n$.
A Poincar\' e series argument shows that $\A^*$ is
isomorphic to $\B^* \otimes F^*$, as a graded $\B^*$-module,
where $F^* = \bigotimes^n_{i=1} \Q[x_i]/(x^{a_i}_i)$. Indeed,
one may extend $\phi$ to a surjective map of graded
$\B^*$--modules,
\begin{equation}
\label{eq=bmod}
\Phi \colon \B^* \otimes F^* \surj \A^*\, ,
\end{equation}
by sending the elements of the monomial $\Q$-basis of $F^*$,
$\{ x_1^{c_1}\cdots x_n^{c_n}\mid 0\le c_i <a_i\, ,\forall i \}$,
to their classes modulo $\I$. By Corollary 2 on p.198 from
\cite{H}, the Poincar\' e series of $\A^*$ is
\[
(1-t^d)^{n-1}\, (1-t^{2kd})\cdot \prod_{i=1}^n (1-t^{w_i})^{-1}\, ,
\]
which is equal to the Poincar\' e series of $\B^*\otimes F^*$.
Therefore, $\Phi$ from \eqref{eq=bmod} above is actually an
isomorphism, which will be used to identify 
$\A^*$ with $\B^*\otimes F^*$.

Set now $V:= \Q- \spn \, \{ y^s_0 \mid 0\le s<k \} \subset \B$,
and $N =V \otimes F$. Clearly, $\dim_{\Q} \A = 2 \dim_{\Q}N$.
The proof ends by checking that $N \cdot N =0$, which is
straightforward.

Part \eqref{nothom}. Assume that $\A^* = H^*(G/K, \Q)$, 
as graded algebras. Write $G/K = \prod^r_{i=1} G_i/K_i$,
with $G_i$ simple, for all $i$. Compute $\der^0$-dimensions
to infer that $G$ must be simple, since $\A$ is a simple $WACI$;
see Proposition~\ref{prop:discr}\eqref{part1} and 
Lemma~\ref{lem=dero}.

By our assumption, $\pi^*_1(\A) = \pi^*_1(H^*(G/K, \Q))$, as
graded vector spaces. One also knows \cite{B} that
$H^*(G/K, \Q)$ is a $WACI$ for which the weights of the variables
are equal to $2$ times the degrees of the fundamental
polynomial invariants of the Weyl group $W_K$, and
likewise for the degrees of the defining relations of
$H^*(G/K, \Q)$.

Case by case inspection of the tables in \cite{Bo} 
reveals that the degree list of the fundamental $W_G$-invariants
contains no element of multiplicity $>2$, if $G$ is simple. For
$n-1>2$, this contradicts the computation from 
Proposition~\ref{prop:discr}\eqref{part3}; see \eqref{eq:min}.
\end{proof}

\subsection{Simple algebras with non-zero signature.}
\label{subsec=nono}
The non-zero signature smoothing problem is more subtle.
For instance, the answer may depend not only on the algebra
structure, but also on the weights, as seen in the following example. 
Let $\A_k$ be $\Q[x]/(x^3)$, where the weight of $x$ is $2k$,
with $k$ odd. It is easy to check that each $\A_k$ is 
a simple $WACI$, with formal dimension $4k$ and signature $\pm 1$.
If $\A_k$ is smoothable, then we infer from 
\cite[Proposition B.5]{Ba} (see also \cite{BLLV}) that
$\sigma (\A_k)$ must be a multiple of $2^{2k-1}-1$. Thus,
$\A_k^*$ is not smoothable, for $k>1$, while
$\A^*_1 = H^*(\C \PP^2, \Q)$ is obviously smoothable.

Let us consider the Eisenbud--Levine $WACI$'s from \cite[p.24]{EL}, with
generators $\{ x_i \}_{1\le i \le n }$ of weight $2$,
and defining relations
\begin{equation}
\label{eq:el}
\begin{cases}
f_1 = x^2_1 - x^2_n \\
\dots \\
f_{n-1} = x^2_{n-1} -x^2_n \\
f_n =x_1 \cdots x_n 
\end{cases}
\end{equation}

Given a $\Q$--$PDA$, $\A$, of formal dimension $4k$, 
we shall denote in the
sequel by $(\A^{2k} ,\cdot_{\omega})\in W(\Q)$ the symmetric inner
product space over $\Q$ obtained from the multiplication map,
$\A^{2k}\otimes \A^{2k} \to \A^{4k}$, via a choice of orientation,
$\omega\in \A^{4k} \setminus \{ 0 \}$.

\begin{prop}
\label{prop:n}

If $\A$ is defined by \eqref{eq:el}, where $n\ge 3$,
then:
\begin{enumerate}
\item \label{p1}
$\A$ is a simple $WACI$.
\item \label{p2}
$\A$ is a $1$-connected Poincar\' e duality algebra,
with formal dimension $4(n-1)$.
\item \label{p3}
There is an orientation, $\omega \in \A^{4(n-1)}\setminus \{ 0 \}$,
such that $(\A^{2(n-1)}, \cdot_{\omega})\in W(\Z)$ and
$\sigma (\A^{2(n-1)}, \cdot_{\omega})= 2^{n-1}$.
\item \label{p4}
The non-trivial odd pseudo-homotopy groups of
$\A$ are:
\[
\begin{cases}
\pi^3_1(\A) =\Q^{n-1}\, , & {\rm and}\\
\pi^{2n-1}_1(\A) = \Q\, . &
\end{cases}
\]
\end{enumerate}
\end{prop}

\begin{proof}
Part \eqref{p1} 
follows from Proposition~\ref{lem:red}, via
\cite[Lemma 2.1(ii)]{PP}. Parts \eqref{p2} and \eqref{p4} are
direct consequences of~\cite[Theorem 3]{H}, and formula
\eqref{eq:min} respectively.

For the proof of Part \eqref{p3}, we will need a convenient
$\Q$-basis of $\A^*$. To obtain it, we first set
$\B :=\Q [y_1, \dots ,y_n]/(y_1-y_n, \dots ,y_{n-1}-y_n) =\Q [y_0]$,
where $\mid y_i \mid =4$, for all $i$, and define a graded algebra
map, 
$\phi \colon \B\to \Q[x_1, \dots, x_n]/(f_1, \dots, f_{n-1})$, by
$\phi (y_i)=x_i^2$, for $i=1, \dots ,n$. 
Like in the proof of Theorem~\ref{thm:smooth}\eqref{mfd}, $\phi$
extends in an obvious way to a surjective map of graded 
$\B^*$--modules,
\begin{equation}
\label{eq=modb}
\Phi \colon \B^*\otimes F^* \surj 
\Q[x_1, \dots, x_n]/(f_1, \dots, f_{n-1})\, ,
\end{equation}
where $F^* =\otimes_{i=1}^n \Q [x_i]/(x_i^2)$.
A standard argument, involving regular sequences and 
Poincar\' e series, shows that the Poincar\' e series of 
$\Q[x_1, \dots, x_n]/(f_1, \dots, f_{n-1})$ is
$(1-t^4)^{n-1}(1-t^2)^{-n}$, which is clearly equal 
to the Poincar\' e series of $\B^*\otimes F^*$.

In this way, we see that the above map \eqref{eq=modb} is an
isomorphism. We will use it to identify $\A^*$ with 
a quotient of $\B^*\otimes F^*$. To be more precise, we will need
the following notation. For $I= \{ i_1<\dots <i_k \}\subset
\{ 1, \dots ,n \}$, we put $x_I :=x_{i_1}\cdots x_{i_k}\in F^{2k}$,
and we denote by $\overline{I}$ the complement of $I$ in
$\{ 1, \dots ,n \}$. It is easy to check in $\B^*\otimes F^*$
the equalities
\[
f_n \cdot x_I \, =\, y_0^{\mid I \mid}\otimes x_{\overline{I}}\, ,
\quad \forall I\, ;
\]
they imply that the images of the monomials
\begin{equation}
\label{eq=monbasis}
\{ y_0^r \otimes x_I \, \mid \, 0\le r< \mid \overline{I} \mid \}
\end{equation}
give a $\Q$-basis of $\A$.
Using this basis, we will show that, if one takes
$\omega =y_0^{n-1}$, then
$(\A, \cdot_{\omega})\in W(\Z)$ and 
$\sigma (\A, \cdot_{\omega})=2^{n-1}$, which clearly proves 
Part~\eqref{p3}, via Poincar\' e duality.

It is not difficult to check that $(\A, \cdot_{\omega})$
splits as an orthogonal direct sum, 
$\bigoplus_{\overline{I}\neq \emptyset} \A_I$, where $\A_I$
denotes the $\Q$-span of the monomials \eqref{eq=monbasis} with
$I$ fixed. For fixed $I$, it is equally easy to check that the
inner product of $y_0^r\otimes x_I$ and $y_0^s\otimes x_I$ equals
\begin{equation}
\label{eq=inner}
\begin{cases}
1\, , & \text{if} \quad r+s=n-1- \mid I \mid \, ;\\
0\, , & \text{otherwise} \, \, .
\end{cases}
\end{equation}

One may use \eqref{eq=inner} to further decompose $\A_I$ into
a sum of hyperbolic planes and one--dimensional summands belonging
to $W(\Z)$, thus checking that $(\A, \cdot_{\omega})\in W(\Z)$. Moreover,
\[
\sigma (\A_I) \, =\, 
\begin{cases}
0\, , & \quad \text{if} \quad n-1- \mid I \mid \quad \text{is odd}\, ;\\
1\, , & \quad \text{if} \quad n-1- \mid I \mid \quad \text{is even}\, ,
\end{cases}
\]
which proves that $\sigma (\A, \cdot_{\omega})= 2^{n-1}$. 
\end{proof}

\begin{remark}
\label{rk=el}
Let $\A$ be a $\Q$--$PDA$, of formal dimension $4k$, with
orientation $\omega \in \A^{4k}\setminus \{ 0 \}$. Then:
the quadratic form $\varphi (\A , \omega)$ has the integrality
property described in \S~\ref{subsec=smpb} if and only if
$(\A^{2k} ,\cdot_{\omega})\in W(\Z)$; see~\cite[Corollary IV.2.6]{MH}.
We should also point out that this integrality result is 
the key fact from our Proposition~\ref{prop:n}\eqref{p3} above
(the value of the signature was computed in~\cite{EL}, using
a different method).
\end{remark}

\begin{theorem}
\label{thm:smooth2}
Let $\A$ be a simple $WACI$ as in 
{\em Proposition~\ref{prop:n}}. Then:
\begin{enumerate}
\item \label{var}
$\A^*$ is the $\Q$-cohomology algebra of a
$1$-connected, smooth closed manifold, $M$. The signature of $M$
equals $\pm 2^{n-1}$.
\item \label{notlie}
$\A^*$ is not isomorphic to $H^*(G/K, \Q)$, for any 
equal rank homogeneous space, $G/K$.
\end{enumerate}
\end{theorem}

\begin{proof}
Part \eqref{var}. We have to show that $\A$ is smoothable.
The integrality condition
is satisfied; see Remark~\ref{rk=el}.

To check the other two $\Q$-surgery obstructions from
Proposition~\ref{prop=surger}, we will construct the total
Pontrjagin class, $q$, as follows. Set $k=n-1$.
Let $p(\C \PP^{2k})= 1+\sum^k_{i=1} c_i u^{2i}$ be
the total Pontrjagin class of $\C \PP^{2k}$, where
$u\in H^2(\C \PP^{2k}, \Q)$ is the canonical generator.
We will take
\begin{equation}
\label{eq:pont}
q := 1+ \sum^k_{i=1} 2^i c_i y_0^i\in \bigoplus_{i\ge 0}
\A^{4i}
\end{equation}
(see \eqref{eq=monbasis}). It readily follows 
from \eqref{eq:pont} that $\A$ has
the same Pontrjagin numbers as $2^k \C \PP^{2k}$. Therefore,
the Hirzebruch signature formula 
\eqref{eq=sign} is also verified for $\A$,
since $\sigma (\A^{2k}, \cdot_{\omega})=2^k$, by
Proposition~\ref{prop:n}\eqref{p3}; see~\cite{M}. 
Hence, $\A$ is smoothable.

Part \eqref{notlie}. The same argument as in the proof of
Theorem~\ref{thm:smooth}\eqref{nothom} gives the result, for
$n\ge 4$. For $n=3$, one may notice that the degree list of
$\pi^*_1(\A)$ is $(3,3,5)$, with $5\not\equiv -1$ (mod $4$),
and resort again to the tables in \cite{Bo}.
\end{proof}

\vspace{0.2cm}

{\sc Stefan Papadima, Inst. of Math. "Simion Stoilow",  P.O. Box 1-764, RO-014700 
Bucharest, Romania;}
{ \it e-mail address:} {\tt Stefan.Papadima@imar.ro}

{\sc Laurentiu Paunescu, School of Maths. and Stats., Univ. of Sydney, Sydney, NSW 2006, 
Australia;} { \it e-mail address:} 
{\tt laurent@maths.usyd.edu.au}

\end{document}